\documentclass{gtmon_a}
\pdfoutput=1

\proceedingstitle{The interaction of finite-type and Gromov--Witten
invariants (BIRS 2003)}
\conferencestart{15 November 2003}
\conferenceend{20 November 2003}
\conferencename{The interaction of finite-type and Gromov--Witten
invariants}
\conferencelocation{Banff International Research Station, Banff, Alberta,
Canada}

\editor{David Auckly}
\givenname{David}
\surname{Auckly}

\editor{Jim Bryan}
\givenname{Jim}
\surname{Bryan}

\title{On the rigidity of stable maps to Calabi--Yau threefolds}

\author{Jim Bryan}
\givenname{Jim}
\surname{Bryan}
\address{Department of Mathematics\\
University of British Columbia\\\newline
Vancouver\\BC\\Canada}
\email{jbryan@math.ubc.ca}
\urladdr{http://www.math.ubc.ca/~jbryan/}

\author{Rahul Pandharipande}
\givenname{Rahul}
\surname{Pandharipande}
\address{Department of Mathematics\\
Princeton University\\\newline
Princeton\\NJ\\USA}
\email{rahulp@math.princeton.edu}
\urladdr{}

\volumenumber{8}
\issuenumber{}
\publicationyear{2006}
\papernumber{4}
\startpage{97}
\endpage{104}

\doi{}
\MR{}
\Zbl{}


\received{19 April 2004}
\revised{2 February 2005}
\accepted{11 March 2005}
\publishedonline{22 April 2006}
\published{22 April 2006}
\proposed{}
\seconded{}
\corresponding{}
\editor{}
\version{}

\arxivreference{math.AG/0405204}

\begin{asciiabstract}
If X is a nonsingular curve in a Calabi--Yau threefold Y whose normal
bundle N_{X/Y} is a generic semistable bundle, are the local
Gromov-Witten invariants of X well defined?  For X of genus two or
higher, the issues are subtle.  We will formulate a precise line of
inquiry and present some results, some positive and some negative.
\\
If $X\subset Y$ is a nonsingular curve in a Calabi--Yau threefold
whose normal bundle $N_{X/Y}$ is a generic semistable bundle, are the
local Gromov--Witten invariants of $X$ well defined? For $X$ of genus
two or higher, the issues are subtle.  We will formulate a precise
line of inquiry and present some results, some positive and some
negative.
\end{asciiabstract}
\AtBeginDocument{\let\bar\wbar\let\tilde\wtilde\let\hat\what}
\def\SU{{\rm SU}}
\makeatletter
\def\cnewtheorem#1[#2]#3{\newtheorem{#1}{#3}[section]
\expandafter\let\csname c@#1\endcsname\c@thm}

\cnewtheorem{theorem}[thm]{Theorem}
\cnewtheorem{lemma}[thm]{Lemma}
\theoremstyle{definition}
\cnewtheorem{definition}[thm]{Definition}
\makeatother  %  move after \newtheorem block

\newcommand{\znums}{{\mathbb Z}}
\newcommand{\End}{\operatorname{End}}
\renewcommand{\P}{\mathbb{P}}
\newcommand{\pic}{\operatorname{Pic}}
\newcommand{\overM}{{}\mskip3mu\overline{\mskip-3mu M\mskip-1mu}\mskip1mu}
\newcommand{\Mh}{\overM_h}
\renewcommand{\O}{\mathcal{O}}

\begin{document}

\begin{htmlabstract}
If X&sub;Y is a nonsingular curve in a Calabi&ndash;Yau threefold
whose normal bundle N<sub>X/Y</sub> is a generic semistable bundle, are the
local Gromov&ndash;Witten invariants of X well defined? For X of genus
two or higher, the issues are subtle.  We will formulate a precise
line of inquiry and present some results, some positive and some
negative.
\end{htmlabstract}

\begin{abstract} 
If $X\subset Y$ is a nonsingular curve in a Calabi--Yau threefold
whose normal bundle $N_{X/Y}$ is a generic semistable bundle, are the
local Gromov--Witten invariants of $X$ well defined? For $X$ of genus
two or higher, the issues are subtle.  We will formulate a precise
line of inquiry and present some results, some positive and some
negative.
\end{abstract}

\maketitle

\section{Introduction}\label{sec-intro}

In 1998, Gopakumar and Vafa \cite{Go-Va-gauge} proposed a duality between $\SU(N)$
Chern--Simons theory on the 3--sphere and topological string theory on the
resolved conifold. As evidence,
Gopakumar and Vafa showed the large--$N$ free energy in
Chern--Simons theory exactly matches (after a change of variables) the
topological string partition function on the resolved conifold.

Mathematically, the {\em topological string partition function} is just the
natural generating function for the Gromov--Witten invariants. The
{\em resolved conifold} is the total space of the bundle $\O _{\P ^{1}}
(-1)\oplus \O _{\P ^{1}} (-1)$, considered as a Calabi--Yau
threefold. The Gromov--Witten theory of the noncompact total space
$\O (-1)\oplus \O (-1)$ 
is well-defined: all
(nonconstant) stable maps have image contained in the zero section
and thus their moduli spaces are compact.

The Gromov--Witten invariants of $\O (-1)\oplus \O (-1)$ are often regarded
as the {\em local Gromov--Witten invariants} of $\P ^{1}$. Indeed, if
$X\subset Y$ is any smoothly embedded rational curve in a Calabi--Yau
threefold $Y$ with normal bundle isomorphic to $\O (-1)\oplus \O
(-1)$, then the contribution of $X$ to the Gromov--Witten invariants of $Y$
is well-defined and is given by the corresponding invariants of $\O
(-1)\oplus \O (-1)$.

We consider here the local theory of higher genus curves. If $X\subset
Y$ is a nonsingular curve in a Calabi--Yau threefold whose normal
bundle $N_{X/Y}$ is a generic semistable bundle, are the local
Gromov--Witten invariants of $X$ well defined? For $X$ of genus two or
higher, the issues are subtle.  We will formulate a precise line of
inquiry and present some results, some positive and some negative.

\subsection*{Acknowledgements}
We thank Michael Thaddeus for his aid in our work.  Jim Bryan is supported
by NSERC, the NSF, and the Sloan foundation; Rahul Pandharipande is
supported by the NSF and the Sloan and Packard foundations.

\section{Definitions and results}

Let $X\subset Y$ be a nonsingular genus--$g$ curve in a threefold
$Y$ with normal bundle $N_{X/Y}$ of
degree $2g{-}2$. If Y is Calabi--Yau, the condition on the normal bundle
is always satisfied.
We define the following notions of rigidity:

\begin{definition}\label{def-rigidity}
\quad 
\begin{itemize}
\item[(i)] A curve $X\subset Y$ is \emph{$(d,h)$--rigid} if for every
degree--$d$, genus--$h$ stable map $f\co C\to X$, we have $H^{0} (C,f^{*}N_{X/Y})=0$. 
\item[(ii)] A curve
$X\subset Y$ is \emph{$d$--rigid} if $X\subset Y$ is
$(d,h)$--rigid for all genera $h$.
\item[(iii)] A curve
$X\subset Y$ is \emph{super-rigid} if $X\subset Y$ is
$d$--rigid for all  $d>0$.
\end{itemize}
\end{definition}

For example, a nonsingular rational curve with normal bundle $\O (-1)\oplus \O
(-1)$ is super-rigid. An elliptic curve $E\subset Y$ is $d$--rigid if and
only if $N_{E/Y}\cong L\oplus L^{-1}$ where $L\to E$ is a flat line bundle
which is not $d$--torsion (see Pandharipande \cite{Pandharipande-degenerate-contributions}). 

For a $(d,h)$--rigid curve $X\subset Y$, the contribution of $X$ to
$N^{h}_{d[X]} (Y)$, the genus--$h$ Gromov--Witten invariant of $Y$ in the
class $d[X]$, is well-defined and given by
\begin{equation}\label{eqn-integrateobstructionbundle}
\int  _{[\Mh (X,d[X])]^\text{vir}}c_{\text{top}} (R^{1}\pi _{*}f^{*}N),
\end{equation}
where $\Mh (X,d)$ is the moduli space of degree--$d$, genus--$h$ stable maps
to $X$,
$$\pi\co U\to \Mh (X,d)$$
is the universal curve,
$$f\co U\to X$$ is the universal map, and $[\ ]^\text{vir}$ denotes
the virtual fundamental class. The $(d,h)$--rigidity of $X$ guarantees
that $R^{1}\pi _{*}f^{*}N$ is a {\em bundle}. See Bryan--Pandharipande
\cite{Br-Pa} for an expanded discussion.

By definition, $(d,h)$--rigidity is a condition on the normal bundle
$N_{X/Y}$.  Assuming $N_{X/Y}$ is generic, we may ask for which pairs
$(d,h)$ does $(d,h)$--rigidity hold.  The $1$--rigidity of a generic normal
bundle is straightforward and was used in
\cite{Pandharipande-degenerate-contributions}.  We prove the following
positive result.

\begin{theorem}\label{thm-genericNisd=Tworigid}
If $X\subset Y$ is a genus--$g$ curve in a threefold $Y$ and $N_{X/Y}$ is a
generic stable bundle of degree $2g{-}2$, then $X$ is $2$--rigid.
\end{theorem}

\noindent However, $3$--rigidity is {\em not} satisfied for genus--$3$ curves.

\begin{theorem}\label{thm-d=Threenotrigid}
If $X\subset Y$ is a genus--3 curve in a threefold $Y$ with
$$\operatorname{deg} (N_{X/Y})=4,$$ then $X$ is not $3$--rigid.
%$(3,h)$--rigid for any $h\geq 10$.
\end{theorem}

Let $N\to X$ be a generic stable bundle of degree $2g{-}2$. By
\fullref{thm-genericNisd=Tworigid}, the degree--2 Gromov--Witten
theory of the total space of $N$ considered as a noncompact threefold is
well-defined by the integral \eqref{eqn-integrateobstructionbundle}.

In the case when $X$ embeds in a threefold $Y$ with normal bundle $N$,
we may regard
the above theory as the degree--2 local Gromov--Witten theory of $X\subset
Y$.  Such embeddings of $X$ can always be found. For example, let $Y$ be
the threefold ${\mathbb P}(\O_X \oplus N)$ with the embedding $X\subset Y$
determined by the trivial factor.  It would be interesting to construct a
curve in a Calabi--Yau threefold with a 2--rigid normal bundle.

The degree--2 {\em local} theory of $X$ is a {\em global} theory for
$N$. Strong global integrality constraints, obtained from the
Gopakumar--Vafa conjecture \cite{Go-Va} and more recently from the
conjectural Gromov--Witten/Donaldson--Thomas correspondence of 
Maulik--Nekrasov--Okounkov--Pandharipande \cite{MNOPOne},
should therefore hold for the degree--2 local theory of $X$.

In \cite{Br-Pa-TQFT,Br-Pa}, a formal local theory is defined in all
degrees for $X\subset Y$. By \fullref{thm-genericNisd=Tworigid}, the above local theory coincides with the degree--2 formal
local theory defined in \cite{Br-Pa-TQFT,Br-Pa}.  However by
\fullref{thm-d=Threenotrigid}, the formal local theory of
\cite{Br-Pa-TQFT,Br-Pa} does {\em not} correspond exactly to a
well-defined global theory of $N$ in degree 3.

Using degeneration techniques, we have recently obtained the complete
computation of the formal local invariants of curves in all degrees
\cite{Br-Pa-local-curves}. The expected integrality holds for the
degree--2 invariants. Integrality for the degree--3 formal local invariants
fails, in fact the breakdown occurs for precisely the {\em same} domain
genus--$h$ for which rigidity fails.

\section[Proof of Theorem \ref{thm-genericNisd=Tworigid}]
{Proof of \fullref{thm-genericNisd=Tworigid}} 

Let $X$ be a nonsingular curve of genus $g$.  Let $N$ be a bundle of rank 2
and degree $2g{-}2$ on $X$.  The bundle $N$ is 2--rigid if
\begin{equation}
\label{eqn-degreeTworigiditycondition}
H^{0} (C,f^{*} (N))=0
\end{equation}
for all stable maps $f\co C\to X$ of degree 2.
If $g=0$, then 
$$N= \O(-1) \oplus \O(-1)$$
is 2--rigid. 
For $g>0$, we will prove 2--rigidity holds 
on an open set of the irreducible 
moduli space of semistable bundles on $X$. As the 
2--rigidity statement was already proven for $g=1$ in 
 \cite{Pandharipande-degenerate-contributions}, we
will assume $g>1$.

Let $\Lambda $ be a fixed line bundle of degree $2g{-}2$. 
Let $\overM_{X} (2,\Lambda )$ be the moduli space of rank--2,
semistable bundles  with determinant $\Lambda $.
Since 2--rigidity is a generic condition, we need only prove
the existence of a 2--rigid bundle $[N] \in \overM_{X} (2,\Lambda )$.

We first prove $H^{0} (X,N)=0$ for an open set
$V\subset \overM_{X} (2,\Lambda ) $. 
Since the vanishing of sections is an open condition
and the moduli space is irreducible, we need only find
an
example. If $L\in \pic^{g-1}(X)$ then the bundle
$$N=L\oplus L^{-1}\Lambda$$
is semistable. Since the locus of $\pic^{g-1}$ determined
by bundles with nontrivial sections is a divisor,
neither $L$ nor $L^{-1}\Lambda $ have sections for generic $L$.
Thus
$$H^0(X, L\oplus L^{-1}\Lambda)=0.$$ 
If there exists a nonzero section $s\in H^0(C,  f^{*} (N))$ for a stable
map $f$, then $s$ must be nonzero on some dominant, irreducible component
of $C'\subset C$.  Also, $s$ must be nonzero when pulled-back to the
normalization of $C'$.  Therefore, to check 2--rigidity for $N$, we
need only prove
$$H^0(C,f^*(N))=0$$
for maps $f\co C\to X$ where the domain is
nonsingular and irreducible. 

Let $C$ be nonsingular and irreducible, and let
$f\co C\rightarrow X$ be a (possibly ramified) double cover. 
Let $B\subset X$ be the branch divisor of $f$. The
map $f$ is well-known to  
determine a square root of $\O (-B) $ in the
Picard group of $X$ by the following construction. 
Let $E=f_{*}
(\O _{C})$. $E$ is a rank--2 bundle on $X$ with a $\znums/2\znums$--action
induced by the Galois group of $f$, and decomposes into a direct sum
$$E\cong \O \oplus Q$$
of $+1$ and $-1$ eigenbundles for the action.
The $+1$ eigenbundle is the trivial line bundle $\O $ while the $-1$
eigenbundle is a line bundle $Q$ with $Q^{2}\cong \O (-B)$.

If $[N] \in V$, a sequence
of isomorphisms
\[
H^{0} (C,f^{*} (N)) \cong
H^{0} (X,f_{*}f^{*} (N)) \cong
H^{0} (X, (\O\oplus Q) \otimes N) \cong
H^{0}
(X,Q\otimes N)
\]
is obtained from the geometry of the 
double cover $f$.
To prove the 2--rigidity of $N$, it therefore suffices to
prove the vanishing of  $H^0(X,Q\otimes N)$ for all 
double covers $f$.

Suppose  $H^{0} (X,Q\otimes N)\neq 0$ for a double cover $f$. 
Then there is a nonzero
sheaf map
$$\iota\co Q^{-1}\rightarrow N.$$ 
Let $S$ be the saturation of the image of $\iota$. 
Then $S$ has the following properties:
\begin{enumerate}
\item[(i)] $S$ is a subbundle of $N$,
%\item[(ii)] $\deg (S)\geq -\deg (Q)=\deg (B)/2$.
\item[(ii)] $S^{2}$ has a section.
\end{enumerate}
Property (ii) is proven as follows. The saturation short exact
sequence of sheaves on $X$ is
\[ 0\to Q^{-1}\to S\to \O _{D}\to 0 \]
for some effective divisor $D$ on $X$. Hence $S\cong Q^{-1} (D)$, and 
$$S^{2}\cong Q^{-2} (2D)\cong \O (B+2D).$$ The line bundle
clearly has a section since both $B$ and $D$ are effective.

For $d\geq 0$,  define $\Delta (d)\subset \overM_{X} (2,\Lambda )$ to
be the following locus in the moduli space:
\[
\bigl\{[N] \ | \ \text{there exists a subbundle}\ S\subset N
\ \text{with}\  \deg (S)=d\ \text{and}\  H^{0} (X,S^{2})\neq 0 \bigr\}
\]
If $[N]\in V$ and $N$ is not 2--rigid, then we have proven  
$[N]$ must lie in $\Delta (d)$ for some $d$.
By \fullref{kskq} below, the proof
 of \fullref{thm-genericNisd=Tworigid} is complete.\qed

\begin{lemma} 
\label{kskq}
For all $d$,
$\dim
\Delta (d)< \dim \overM_{X} (2,\Lambda ).$
\end{lemma}

\begin{proof}
By stability, $\Delta (d)$ is empty if $d\geq g-1$. We may assume 
\[ 0\leq d< g-1.  \]
If $[N]\in \Delta (d)$, then $N$ is given by an extension
\[ 0\to S\to N\to S^{-1}\otimes \Lambda \to 0. \]
The number of parameters for $S$ is $g$, the dimension of
$\pic^d (X)$. The number of
parameters of the space of extensions is 
\[
h^{1} (X,S^{2}\otimes \Lambda ^{-1})-1=3g-3-2d-1,
\]
by Riemann--Roch. Therefore
\[
\dim \Delta (d)\leq g+ 3g-3-2d-1=4g-4-2d.
\]
If $d> (g-1)/2$, then $\dim \Delta (d)<\dim M_{X} (2,\Lambda )=3g-3$.

By the above analysis, we may assume $0\leq d\leq (g-1)/2$.  We will
now redo the dimension count for $\Delta(d)$ using the condition
$H^{0} (S^{2})\neq 0$. By Brill--Noether theory, the dimension of the
space of line bundles of degree $2d<g$ having nonzero sections is
$2d$. Recomputing, we find
\[
\dim \Delta (d)\leq 2d + 3g-3 -2d-1 = 3g-4.
\]
Therefore $\dim \Delta (d)<\dim \overM_{X} (2,\Lambda )$ for all $d$.
\end{proof}

\section[{Proof of Theorem \ref{thm-d=Threenotrigid}}]
{Proof of \fullref{thm-d=Threenotrigid}}

Let $X$ be a nonsingular genus--3 curve.  For any rank--2 bundle 
$N$ on $X$ of degree 4, we will show there exists a degree--3 stable map
$$f\co C\to X$$
for which $H^{0} (C,f^{*} (N))\neq 0$.

If $N$ contains a rank--1 subbundle $L\subset N$ of positive
degree,  then $H^{0} (X,L^{3})\neq 0$ by Riemann--Roch. Let $s$ be a nonzero
section of $L^{3}$ and let $f\co C\to X$ 
be the cube root of $s$ in $L$. Since
$H^{0} (C,f^{*} (L)) $ has a canonical section, we have shown
 $H^{0}
(C,f^{*} (N))\neq 0$.
If $N$ is unstable, then it has a destabilizing 
rank--1 subbundle of positive degree, and is not 3--rigid.

Let $\Lambda$ be a line bundle on $X$ of degree 4.
Let $\overM_X(2, \Lambda)$ be the moduli space of rank--2, semistable
bundles with determinant $\Lambda$.
We will prove, for general $$[N]\in \overM_X(2, \Lambda),$$ 
that $N$ contains a rank--1 subbundle of degree 1.
Hence the general semistable bundle $N$ with determinant $\Lambda$
is not 3--rigid.
By semicontinuity and the irreducibility of the moduli
of semistable bundles, we conclude no semistable bundle
on $X$ is 3--rigid.

The dimension of the moduli space
$\overM_{X} (2,\Lambda )$ is 6. The extensions 
\begin{equation}\label{eqn-exactseqforN}
0\to L\to N\to L^{-1}\Lambda \to 0,
\end{equation}
where $L$ has rank 1 and degree 1, are parametrized (up to scale) by a
projective bundle $B$ with fiber  
$\P (H^{1} (X,L^{2}\Lambda^{-1}))$ over $\pic^1(X)$.
Since the dimensions of $\pic^1(X)$ and $\P (H^{1}
(X,L^{2}\Lambda^{-1} ))$ are both 3, the dimension of $B$ is 6. 
An elementary argument shows the map to moduli,
\[
\epsilon\co  B  \rightarrow \overM_{X} (2,\Lambda),
\]
is well-defined on a (nonempty) open set of $B$. If $\epsilon$
is dominant, the theorem is proven.

It suffices to prove that
the tangent space to $B$ generically
surjects onto the tangent space of $\overM_{X} (2,\Lambda )$.
The following
argument was provided by Michael Thaddeus.

Let $\End _{0}N$ be the bundle of traceless endomorphisms of $N$. The
tangent space to $\overM_{X} (2,\Lambda )$ at $N$ is given by
$H^{1} (X,\End _{0}N)$. The exact sequence \eqref{eqn-exactseqforN} induces a filtration
\[
0=E_{0}\subset E_{1}\subset E_{2}\subset E_{3}=\End _{0}N
\]
where 
\[
E_{1}/E_{0}=L^{2}\Lambda ^{-1},\quad E_{2}/E_{1}=\O ,\quad\text{and}
\ E_{3}/E_{2}=L^{-2}\Lambda .
\]
For generic $L$, 
\[
H^{0} (X,L^{-2}\Lambda )=H^{1} (X,L^{-2}\Lambda )=0.
\]
Hence,
the inclusion $E_{2}\subset E_{3}$ induces isomorphisms
\[
H^{1} (X,E_{2})\cong H^{1} (X,\End _{0}N) \quad\text{and}\quad H^{0} (X,E_{2})\cong H^{0}
(X,\End _{0}N)=0.
\]
Consequently, the exact sequence for the inclusion $E_{1}\subset E_{2}$
reduces to
\[
0\to H^{1} (X,L^{2}\Lambda ^{-1})/H^{0} (X,\O )\to H^{1} (X,\End _{0}N)\to
H^{1} (X,\O )\to 0.
\]
The tangent space of $\overM_{X} (2,\Lambda )$ at $N$
is therefore identified
with the tangent space of $B$ at $N$. Indeed, $H^{1} (X,L^{2}\Lambda
^{-1})/H^{0} (X,\O )$ is the space of deformations of the extension class
and $H^{1} (X,\O )$ is the space of deformations of $L$.  \qed

\bibliographystyle{gtart}
\bibliography{link} 

\end{document}